\documentclass[12pt,leqno]{amsart}
\usepackage{amsmath}
\usepackage{amssymb}

\def\underset#1#2{{\mathrel{\mathop {{}_{} {#2}}\limits_{{#1}_{}}}}}
\def\upplim_#1{\underset{#1}{\overline\lim}\;}
\def\lowlim_#1{\underset{#1}{\underline\lim}\;}
\newtheorem{corollary}[equation]{Corollary}
\newtheorem{definition}[equation]{\indent{\it Definition}\rm }
\newtheorem{claim}[equation]{\indent{\it Claim}\rm }

\newtheorem{lem}[equation]{Lemma}
\newtheorem{lemma}[equation]{Lemma}
\newtheorem{proposition}[equation]{Proposition}

\newtheorem{theorem}[equation]{Theorem}

\newcommand{\C}{{\mathbf{C}}}

\renewcommand{\P}{{\mathbf{P}}}

\newcommand{\R}{{\mathbf{R}}}

\newcommand{\rank}{\mathrm{rank}}
\newcommand{\zero}{\mathrm{Zero}}

\newcommand{\supp}{\mathrm{Supp}\,}

\newcommand{\Z}{\mathbf{Z}}

\numberwithin{equation}{section}

\title[Second main theorems for meromorphic mappings]{Second main theorems for meromorphic mappings intersecting moving hyperplanes with truncated counting functions and unicity problem} 

\date { }
\author{Si Duc Quang}

\address{Department of Mathematics, Hanoi National University of Education\\
136-Xuan Thuy, Cau Giay, Hanoi, Vietnam}
\email{Email: quangsd@hnue.edu.vn}

\begin{document}

\begin{abstract}
In this article, we show some new second main theorems for the mappings and moving hyperplanes of $\P^n(\C)$ with truncated counting functions. Our results are improvements of recent previous second main theorems for moving hyperplanes with the truncated (to level $n$) counting functions. 
%In particular, we show that $$\parallel\dfrac{2q}{3(n+1)}T_f(r) \le \sum_{i=1}^q N_{(f,a_i)}^{[n]}(r) + o(T_f(r)) + O(\max_{1\le i \le q}T_{a_i}(r))\text{ if $q\ge 3n+3$},$$$$\text{and }\parallel \dfrac{q-n+1}{n+2}T_f(r) \le \sum_{i=1}^q N_{(f,a_i)}^{[n]}(r) + o(T_f(r)) + O(\max_{1\le i \le q}T_{a_i}(r))\text{ if $2n+1\le q\le 3n+2.$}$$
As their application, we prove a unicity theorem for meromorphic mappings sharing moving hyperplanes.
\end{abstract}

\def\thefootnote{\empty}
\footnotetext{
2010 Mathematics Subject Classification:
Primary 32H30, 32A22; Secondary 30D35.\\
\hskip8pt Key words and phrases: Nevanlinna, second main theorem, meromorphic mapping, moving hyperplane.}

\maketitle

\section{Introduction}

\noindent
The theory of the Nevanlinna's second main theorem for meromorphic mappings of $\C^m$ into the complex projective space $\P^n(\C)$ intersecting a finite set of fixed hyperplanes or moving hyperplanes in $\P^n(\C)$ was started about 70 years ago and  has grown into a huge theory. For the case of fixed hyperplanes, maybe, the second main theorem given by Cartan-Nochka is the best possible. Unfortunately, so far there has been a few second main theorems with truncated counting functions for moving hyperplanes. Moreover, almost of them are not sharp.

We state here some recent results on the second main theorems for moving hyperplanes with truncated counting functions. 
 
Let $\{a_i\}_{i=1}^q $ be  meromorphic mappings of $\C^m$ into the dual space $\P^n(\C)^*$ in general position. For the case of nondegenerate meromorphic mappings, the second main theorem with truncated (to level $n$) counting functions states that.

\noindent
\textbf{Theorem A} (see \cite[Theorem 2.3]{MR} and \cite[Theorem 3.1]{TQ05}). {\it Let $f :\C^m \to \P^n(\C)$ be a meromorphic mapping. Let $\{a_i\}_{i=1}^q \ (q\ge n+2)$ be  meromorphic mappings of $\C^m$ into $\P^n(\C)^*$ in general position such that $f$ is linearly nondegenerate over $\mathcal {R}(\{a_i\}_{i=1}^{q}).$ Then
$$|| \ \dfrac {q}{n+2}  T_f(r) \le \sum_{i=1}^q N_{(f,a_i)}^{[n]}(r) + o(T_f(r)) + O(\max_{1\le i \le q}T_{a_i}(r)).$$ }

We note that, Theorem A is still the best second main theorem with truncated counting functions for nondegenerate meromorphic mappings and moving hyperplanes available at present. In the case of degenerate meromorphic mappings, the second main theorem for moving hyperplanes with counting function truncated to level $n$ was first given by M. Ru-J. Wang \cite{RW} in 2004. After that in 2008, D. D. Thai-S. D. Quang \cite{TQ08} improved the result of M. Ru-J. Wang by proved the following second main theorem.

\noindent
\textbf{Theorem B} (see \cite[Corollary 1]{TQ08}). {\it Let $f:\C^m\to\P^n(\C)$ be a meromorphic mapping. Let $\{a_i\}_{i=1}^q$ $(q\ge 2n+1)$ be $q$ meromorphic mappings of $\C^m$ into $\P^n(\C)^*$ in general position such that $(f,a_i)\not\equiv 0\ (1\le i\le {q}).$  
Then
 $$\bigl |\bigl |\quad\dfrac{q}{2n+1}\cdot T_f(r)\le \sum_{i=1}^{q}N^{[n]}_{(f,a_i)}(r)+O\bigl(\max_{1\le i\le {q}}T_{a_i}(r)\bigl)+O\bigl(\log^{+}T_f(r)\bigl).$$ }

These results play very essential roles in almost all researches on truncated multiplicity problems of  meromorphic mappings with moving hyperplanes. Hovewer, in our opinion, the above mentioned results of these authors are still weak.

Our main purpose of the present paper is to show a stronger second main theorem of  meromorphic mappings from $\C^m$ into $\P^n(\C)$ for moving targets. Namely, we will prove the following.
\begin{theorem}\label{1.1}
Let $f :\C^m \to \P^n(\C)$ be a meromorphic mapping. Let $\{a_i\}_{i=1}^q \ (q\ge 2n-k+2)$ be meromorphic mappings of $\C^m$ into $\P^n(\C)^*$ in general position such that $(f,a_i)\not\equiv 0\ (1\le i\le q),$ where $k+1=\rank_{\mathcal R\{a_i\}}(f)$. Then the following assertions hold:
\begin{align*}
&\mathrm{(a)}\ || \ \dfrac {q}{2n-k+2}T_f(r) \le \sum_{i=1}^q N_{(f,a_i)}^{[k]}(r) + o(T_f(r)) + O(\max_{1\le i \le q}T_{a_i}(r)),\\
&\mathrm{(b)}\ || \ \dfrac {q-n+2k-1}{n+k+1}T_f(r) \le \sum_{i=1}^q N_{(f,a_i)}^{[k]}(r) + o(T_f(r)) + O(\max_{1\le i \le q}T_{a_i}(r)).
\end{align*}
\end{theorem}
We may see that Theorem \ref{1.1}(a) is a generalization of Theorem A and also is an  improvement of Theorem B. Theorem \ref{1.1}(b) is really stronger than Theorem B. 

\noindent
\textit{Remark.} %1) The following example shows that, for any $2n-k+1\ (1\le k\le n)$ hyperplanes $H_1,\ldots ,H_{2n-k+1}$ of $\P^n(\C)$ in general position, there exist a nonconstant holomorphic curve $f:\C\rightarrow \P^n(\C)$ with $\rank_{\C}(f)=k+1$ such that $f$ does not intersect these all hyperplanes. Therefore, the condition ``$q\ge 2n-k+2$''  in Theorem \ref{1.1} can not be removed.

%Indeed, let $O=\bigcap_{i=1}^nH_i$ and $H=\bigcap_{i=n+1}^{2n-k+1}H_i$. Since $H_1,\ldots ,H_{2n-k+1}$ are in general position, we see that $O$ is a point, $H$ is a linear subspace of $\P^n(\C)$ of dimiension $k-1$ and $O\not\in H$. Take $L$ the smallest linear subspace of $\P^n(\C)$ containing $O$ and $H$. Then we may regard $L$ as a projective space of dimension $k$ and $H$ as a hyperplane of $L$.

%It is enough to show that there exists a linealy non-degenerate holomorphic curve $f:\C\rightarrow L$ omitting $O$ and $H$. Chose a homogeneous coordinates $(\omega_0:\cdots \omega_k)$ so that $O$ has coordinate $(1:0:\cdots :0)$ and $H$ is given by
%$$ H=\{(\omega_0:\cdots :\omega_k);\ \omega_0=0\}.$$ 
%We consider the holomorphic curve $f$ from $\C$ into $L$ defined by
%$$ f(z)=(1:e^z:e^{z^2}:\cdots :e^{z^k}). $$
%It is clear that $f$ omits the point $O$ and the hyperplane $H$. We also see that $f$ is linearly non-degenerate.

1) If $k\ge\dfrac{n+1}{2}$ then Theorem \ref{1.1}(a) is stronger than Theorem \ref{1.1}(b). Otherwise, if $k<\dfrac{n+1}{2}$ then Theorem \ref{1.1}(b) is stronger than Theorem \ref{1.1}(a).

2) If $k=0$ then $f$ is constant map, and hence $T_f(r)=0.$

3) Setting $t=\frac{2n-k+2}{3n+3}$ and $\lambda =\frac{n+k+1}{3n+3},$ we have $t+\lambda =1$. Thus, for all $1\le k\le n$ we have
\begin{align*}
\max\biggl \{\dfrac{q}{2n-k+2},\dfrac{q-n+2k-1}{n+k+1}\biggl\}&\ge \dfrac{q}{2n-k+2}\cdot t+\dfrac{q-n+2k-1}{n+k+1}\cdot\lambda\\
\\
&=\dfrac{2q-n+2k-1}{3n+3}\ge\dfrac{2q-n+1}{3n+3}.
\end{align*}

4) If $k\ge 1$, we have the following estimates:
\begin{itemize}
\item $\min_{\frac{n+1}{2}\le k\le n, (k\in\Z)}\left (\dfrac{q}{2n-k+2}\right )\ge\dfrac{q}{2n-\frac{n+1}{2}+2}=\dfrac{2q}{3(n+1)}$.
\item $\min_{1\le k\le \frac{n+1}{2}, (k\in\Z)}\left (\dfrac{q-n+2k-1}{n+k+1}\right )=\min_{1\le k\le \frac{n+1}{2}, (k\in\Z)}\left (\dfrac{q-3n-3}{n+k+1}+2\right )$

\noindent
\hspace{185pt}$\ge\begin{cases}
\dfrac{2q}{3(n+1)}&\text{ if }q\ge 3n+3 \\
\dfrac{q-n+1}{n+2}&\text{ if }q< 3n+3
\end{cases}$
\end{itemize}
Thus
$$\min_{1\le k\le n}\biggl \{\max\bigl \{\dfrac{q}{2n-k+2},\dfrac{q-n+2k-1}{n+k+1}\bigl\}\biggl \}
\ge\begin{cases}
\dfrac{2q}{3(n+1)}&\text{ if }q\ge 3n+3 \\
\dfrac{q-n+1}{n+2}&\text{ if }q< 3n+3.
\end{cases}$$
Therefore, from Theorem \ref{1.1} and Remark (1-4) we have the following corollary.
\begin{corollary}\label{1.2}
Let $f :\C^m \to \P^n(\C)$ be a meromorphic mapping. Let $\{a_i\}_{i=1}^q \ (q\ge 2n+1)$ be  meromorphic mappings of $\C^m$ into $\P^n(\C)^*$ in general position such that $(f,a_i)\not\equiv 0\ (1\le i\le q).$ 

$\mathrm{(a)}$ Then we have
$$|| \dfrac{2q-n+1}{3(n+1)}T_f(r) \le \sum_{i=1}^q N_{(f,a_i)}^{[n]}(r) + o(T_f(r)) + O(\max_{1\le i \le q}T_{a_i}(r)).$$

$\mathrm{(b)}$ If $q\ge 3n+3$ then 
$$|| \dfrac{2q}{3(n+1)}T_f(r) \le \sum_{i=1}^q N_{(f,a_i)}^{[n]}(r) + o(T_f(r)) + O(\max_{1\le i \le q}T_{a_i}(r)).$$

$\mathrm{(c)}$ If $q< 3n+3$ then 
$$|| \dfrac{q-n+1}{n+2}T_f(r) \le \sum_{i=1}^q N_{(f,a_i)}^{[n]}(r) + o(T_f(r)) + O(\max_{1\le i \le q}T_{a_i}(r)).$$
\end{corollary}

As applications of these second main theorems, in the last section we will prove a unicity theorem for meromorphic mappings sharing moving hyperplanes regardless of multiplicities. To state our main result, we give the following definition.
 
Let $f:\C^m \to \P^n(\C)$ be a meromorphic mapping. Let $k$ be a positive integer or maybe $+\infty$. Let $\{a_i\}_{i=1}^{q}$ be ``slowly'' (with respect to $f$) moving hyperplanes in $\P^n(\C)$ in general position such that 
$$\dim\ \{z\in\C^m: (f,a_i)(z)\cdot (f,a_j)(z)= 0\} \leq m-2\quad (1\le i<j\le q).$$

Consider the set $\mathcal F(f,\{a_i\}_{i=1}^q,k)$ of all meromorphic maps $g: \C^m \to \P^n(\C)$ satisfying the following two conditions:

(a) $\min\{\nu_{(f,a_i)}(z),k\}=\min\{\nu_{(g,a_i)}(z),k\}\quad (1 \le i \le q),$ for all $z\in\C^m$,

(b) $f(z) = g(z)$ for all $z\in\bigcup_{i=1}^{q}\zero (f,a_i)$.

We wil prove the following

\begin{theorem}\label{1.3}
Let $f:\C^m \to \P^n(\C)$ be a meromorphic mapping. Let $\{a_i\}_{i=1}^{q}$ be slowly (with respect to $f$) moving hyperplanes in $\P^n(\C)$ in general position such that 
$$\dim\ \{z\in\C^m: (f,a_i)(z)\cdot (f,a_j)(z)= 0\} \leq m-2\quad (1\le i<j\le q).$$
Then the following assertions hold:

a) If $q>\frac{9n^2+9n+4}{4}$ then $\sharp\ \mathcal F(f,\{a_i\}_{i=1}^q,1)\le 2,$ 

b) If $q>3n^2+n+2$ then $\sharp\ \mathcal F(f,\{a_i\}_{i=1}^q,1)=1.$
\end{theorem}

{\bf Acknowledgements.} This work was done during a stay of the author at Vietnam Institute for Advanced Study in Mathematics. He would like to thank the institute for their support.

\section{Basic notions and auxiliary results from Nevanlinna theory}

\noindent
{\bf (a)} Counting function of divisor.

For $z = (z_1,\dots,z_m) \in \C^m$, we set
$\Vert z \Vert = \Big(\sum\limits_{j=1}^m |z_j|^2\Big)^{1/2}$
and define 
\begin{align*}
B(r) &= \{ z \in \C^m ; \Vert z \Vert < r\},\quad
S(r) = \{ z \in \C^m ; \Vert z \Vert = r\},\\
d^c & = \dfrac{\sqrt{-1}}{4\pi}(\overline \partial - \partial),\quad
\sigma = \big(dd^c \Vert z\Vert^2\big)^{m-1},\\
\eta &= d^c \text{log}\Vert z\Vert^2 \land 
\big(dd^c\text{log}\Vert z \Vert\big)^{m-1}.
\end{align*}

Thoughout this paper, we denote by $\mathcal M$ the set of all meromorphic functions on $\C^m$. A divisor $E$ on $\C^m$ is given by a formal sum $E=\sum\mu_{\nu}X_{\nu}$, 
where $\{X_\nu\}$ is a locally family of distinct irreducible analytic hypersurfaces in $\C^m$ and $\mu_{\nu}\in \Z$. We define the support of the divisor
$E$ by setting $\supp (E)=\cup_{\nu\ne 0} X_\nu$.
Sometimes, we identify the divisor $E$ with a function $E(z)$ from $\C^m$ 
into $\Z$ defined by $E(z):=\sum_{X_{\nu}\ni z}\mu_\nu$. 

Let $k$ be a positive integer or $+\infty$. We define the truncated divisor $E^{[k]}$ by 
$$
E^{[k]}:= \sum_{\nu}\min\{\mu_\nu, k \}X_\nu ,
$$ 
and the {\it truncated counting function to level $k$} of $E$ by
\begin{align*}
N^{[k]}(r,E) := \int\limits_1^r \frac{n^{[k]}(t,E)}{t^{2m-1}}dt\quad 
(1 < r < +\infty),
\end{align*}
where 
\begin{align*}
n^{[k]}(t,E): =
\begin{cases}
\int\limits_{\supp (E) \cap B(t)} E^{[k]}\sigma &\text{ if } m \geq 2,\\
\sum_{|z| \le t} E^{[k]}(z)&\text{ if } m = 1.
\end{cases}
\end{align*}
We omit the character $^{[k]}$ if $k=+\infty$.

For an analytic hypersurface $E$ of $\C^m$, we may consider it as a reduced divisor and denote by $N(r,E)$ its counting function.

Let $\varphi$ be a nonzero meromorphic function on $\C^m$.
We denote by $\nu^0_{\varphi}$ (resp. $\nu^{\infty}_{\varphi}$) the divisor of zeros (resp. divisor of poles) of $\varphi$. The divisor of $\varphi$ is defined by
$$\nu_{\varphi}=\nu^0_{\varphi}-\nu^{\infty}_{\varphi}.$$

We have the following Jensen's formula:
\begin{align*}
N(r,\nu^0_{\varphi}) - N(r,\nu^{\infty}_{\varphi}) =
\int\limits_{S(r)} \text{log}|\varphi| \eta 
- \int\limits_{S(1)} \text{log}|\varphi| \eta .
\end{align*}
For convenience, we will write $N_{\varphi}(r)$ and $N^{[k]}_{\varphi}(r)$ for $N(r,\nu^0_{\varphi})$ and $N^{[k]}(r,\nu^0_{\varphi})$, respectively.

\noindent
{\bf (b)} The first main theorem.

Let $f$ be a meromorphic mapping of $\C^m$ into $\P^n(\C)$. For arbitrary fixed homogeneous coordinates $(w_0: \cdots : w_n)$ of $\P^n(\C )$, we take a reduced representation $f = (f_0 : \cdots : f_n)$, which means that each $f_i$ is holomorphic function on $\C^m$ and $f(z) = (f_0(z) : \cdots : f_n(z))$ outside the analytic set $I(f):=\{ z ; f_0(z) = \cdots = f_n(z) = 0\}$ of codimension at least $2$.

Denote by $\Omega$ the Fubini Study form of $\P^n(\C )$. The characteristic function of $f$ (with respect to $\Omega$) is defined by
\begin{align*}
T_f(r) := \int_1^r\dfrac{dt}{t^{2m-1}}\int_{B(t)}f^*\Omega\wedge\sigma ,\quad\quad 1 < r < +\infty.
\end{align*}
By Jensen's formula we have
\begin{align*}
T_f(r)=\int_{S(r)}\log ||f||\eta +O(1),
\end{align*}
where $\Vert f  \Vert = \max \{ |f_0|,\dots,|f_n|\}$.

Let $a$ be a meromorphic mapping of $\C^m$ into $\P^n(\C)^*$ with reduced representation
$a = (a_0 : \dots : a_n)$. We define
$$m_{f,a}(r)=\int\limits_{S(r)} \text{log}\dfrac {||f||\cdot ||a||}{|(f,a)|}\eta -
\int\limits_{S(1)}\text{log}\dfrac {||f||\cdot ||a||}{|(f,a)|}\eta,$$
where $\Vert a \Vert = \big(|a_0|^2 + \dots + |a_n|^2\big)^{1/2}$ and $(f,a)=\sum_{i=0}^nf_i\cdot a_i.$

Let $f$ and $a$ be as above. If $(f,a)\not \equiv 0$, then the first main theorem for moving hyperplaness in value distribution theory states
$$T_f(r)+T_a(r)=m_{f,a}(r)+N_{(f,a)}(r)+O(1)\ (r>1).$$

For a meromorphic function $\varphi$ on $\C^m$, the proximity function $m(r,\varphi)$ is defined by
$$ m(r,\varphi) = \int\limits_{S(r)} \log^+ |\varphi| \eta , $$
where $\log^+ x = \max \big\{ \log x, 0\big\}$ for $x \geqslant 0$.
The Nevanlinna's characteristic function is defined by
$$T(r, \varphi ) = N(r, \nu^{\infty}_\varphi) + m(r,\varphi ).$$
We regard $\varphi$ as a meromorphic mapping of $\C^m$ into $\P^1(\C )^*$, there is a fact that
$$ T_\varphi (r)=T(r,\varphi )+O (1). $$
{\bf (c)} Lemma on logarithmic derivative.

 As usual, by the notation $``|| \ P"$  we mean the assertion $P$ holds for all $r \in [0,\infty)$ excluding a Borel subset $E$ of the interval $[0,\infty)$ with $\int_E dr<\infty$. Denote by $\Z_+$ the set of all nonnegative integers. The lemma on logarithmic derivative in Nevanlinna theorey is stated as follows.
\begin{lem}[{see \cite[Lemma 3.11]{Shi}}]\label{2.1}
Let $f$ be a nonzero meromorphic function on $\C^m.$ Then 
$$\biggl|\biggl|\quad m\biggl(r,\dfrac{\mathcal{D}^\alpha (f)}{f}\biggl)=O(\log^+T_f(r))\ (\alpha\in \Z^m_+).$$
\end{lem}

\noindent
{\bf (d)} Family of moving hyperplanes. 

We assume that thoughout this paper, the homogeneous coordinates of $\P^n(\C)$ is chosen so that for each given meromorphic mapping $a=(a_0:\cdots :a_n)$ of $\C^m$ into $\P^n(\C)^*$ then $a_{0}\not\equiv 0$. We set
$$ \tilde a_i=\dfrac{a_i}{a_0}\text{ and }\tilde a=(\tilde a_0:\tilde a_1:\cdots:\tilde a_n).$$
Let $f:\C^m\rightarrow\P^n(\C)$ be a meromorphic mapping with the reduced representation $f=(f_0:\cdots :f_n).$  We put $(f,a):=\sum_{i=0}^{n}f_ia_{i}$ and $(f,\tilde a):=\sum_{i=0}^{n}f_i\tilde a_{i}.$

Let $\{a_i\}_{i=1}^q$ be $q$ meromorphic mappings of $\C^m$ into $\P^n(\C)^*$  with reduced representations $a_i=(a_{i0}:\cdots :a_{in})\ (1\le i\le q).$ We denote by  $\mathcal R(\{a_i\})$ (for brevity we will write $\mathcal R$ if there is no confusion) the smallest subfield of $\mathcal M$ which contains $\C$ and all ${a_{i_j}}/{a_{i_k}}$ with $a_{i_k}\not\equiv 0.$

%Let $N\ (N\ge n)$ be an integer. We say that the family $\{a_i\}_{i=1}^q$ is in \textit{$N-$subgeneral position}  if $\dim (\{a_{i_0},\ldots ,a_{i_N}\})_{\mathcal {M}}=n+1$ for any $1\le i_0\le\cdots\le i_N$, where $(\{a_{i_0},\ldots ,a_{i_N}\})_{\mathcal {M}}$ is the linear span of $\{a_{i_0},\ldots ,a_{i_N}\}$ over the field $\mathcal{ M}.$

%If $\{a_i\}_{i=1}^q$ is in $n-$subgeneral position then we say that it is in \textit{general position}

\begin{definition}
The family $\{a_i\}_{i=1}^q$ is said to be in general position  if $\dim (\{a_{i_0},\ldots ,a_{i_n}\})_{\mathcal {M}}=n+1$ for any $1\le i_0\le\cdots\le i_n\le q$, where $(\{a_{i_0},\ldots ,a_{i_n}\})_{\mathcal {M}}$ is the linear span of $\{a_{i_0},\ldots ,a_{i_N}\}$ over the field $\mathcal{ M}.$
\end{definition}

\begin{definition}
A subset $\mathcal {L}$ of $\mathcal {M}$ (or $\mathcal {M}^{n+1}$) is said to be   minimal over the field $\mathcal R$ if it is linearly dependent over $\mathcal {R}$ and each proper subset of $\mathcal L$ is linearly independent over  $\mathcal {R}.$
\end{definition}

Repeating the argument in (\cite[Proposition 4.5]{Fu}), we have the following:
\begin{proposition}[{see \cite[Proposition 4.5]{Fu}}]\label{2.2}
Let $\Phi_0,\ldots,\Phi_k$ be meromorphic functions on $\C^m$ such that $\{\Phi_0,\ldots,\Phi_k\}$ 
are  linearly independent over $\C.$
Then  there exists an admissible set $\{\alpha_i=(\alpha_{i1},\ldots,\alpha_{im})\}_{i=0}^k \subset \Z^m_+$ with $|\alpha_i|=\sum_{j=1}^{n}|\alpha_{ij}|\le k \ (0\le i \le k)$ such that the following are satisfied:

(i)\  $\{{\mathcal D}^{\alpha_i}\Phi_0,\ldots,{\mathcal D}^{\alpha_i}\Phi_k\}_{i=0}^{k}$ is linearly independent over $\mathcal M,$\ i.e, \ 
$\det{({\mathcal D}^{\alpha_i}\Phi_j)}\not\equiv 0.$ 

(ii) $\det \bigl({\mathcal D}^{\alpha_i}(h\Phi_j)\bigl)=h^{k+1}\det \bigl({\mathcal D}^{\alpha_i}\Phi_j\bigl)$ for
any nonzero meromorphic function $h$ on $\C^m.$
\end{proposition}

\section{Proof of Theorem \ref{1.1}}

In order to prove Theorem \ref{1.1} we need the following.
\begin{lem}\label{3.1}
Let  $f:\C^m\rightarrow\P^n(\C)$ be a meromorphic mapping. Let $\{a_i\}_{i=1}^q$ $(q\ge n+1)$ be $q$ meromorphic mappings of $\C^m$ into $\P^n(\C)^*$ in general position. Assume that there exists a partition $\{1,\ldots,q\}=I_1\cup I_2\cdots\cup I_l$ satisfying:

$\mathrm{(i)}$ \ $\{(f,\tilde a_i)\}_{i\in I_1}$ is minimal over $\mathcal R$, and $\{(f,\tilde a_i)\}_{i\in I_t}$ is linearly independent over $\mathcal {R}\ (2\le t \le l), $ 

$\mathrm{(ii)}$ \ For any $2\le t\le l,i\in I_t,$ there exist meromorphic functions $c_i\in \mathcal {R}\setminus\{0\}$ such that 
$$\sum_{i\in I_t}c_i(f,\tilde a_i)\in \biggl(\bigcup_{j=1}^{t-1}\bigcup_{i\in I_j}(f,\tilde a_i) \biggl)_{\mathcal {R}}.$$
Then we have
$$ T_f(r)\le\sum_{i=1}^qN^{[k]}_{(f,a_i)}+ o(T_f(r)) + O(\max_{1\le i \le q}T_{a_i}(r)),$$
where $k+1=\rank_{\mathcal R}(f)$.
\end{lem}
\textbf{Proof.}\ Let $f=(f_0:\cdots :f_n)$ be a reduced representation of $f$. By changing the homogeneous coordinate system of $\P^n(\C)$ if necessary, we may assume that $f_0\not\equiv 0.$ 
Without loss of generality, we may assume that $I_1=\{1,\ldots.,k_1\}$ and
$$I_t=\{k_{t-1}+1,\ldots, k_t\}\ (2\le t \le l),\text{ where }1=k_0<\cdots< k_l=q.$$ 

Since $\{(f,\tilde a_i)\}_{i\in I_1}$ is minimal over $\mathcal R$, there exist $c_{1i}\in\mathcal {R}\setminus \{0\}$ such that 
$$\sum_{i=1}^{k_1}c_{1i}\cdot (f, \tilde a_i)=0.$$
Define $c_{1i}=0$ for all $i>k_1.$ Then 
$$\sum_{i=1}^{k_l} c_{1i}\cdot (f,\tilde a_i)=0.$$
Because ${\{c_{1i}(f,\tilde a_i)\}}_{i=k_0+1}^{k_1}$ is linearly independent over $\mathcal R,$ Lemma \ref{2.2} yields that there exists an admissible set $\{\alpha_{1(k_0+1)},\ldots,\alpha_{1k_1}\}\subset \Z^m_+$ \ $(|\alpha_{1i}|\le k_1-k_0-1\le \rank_{\mathcal R}f-1=k)$ such that the matrix
$$\ A_1=\left (\mathcal {D}^{\alpha_{1i}}(c_{1j}(f,\tilde a_j));k_0+1\le i,j\le k_1 \right)$$
has nonzero determinant.

Now consider $t\ge 2.$ 
By constructing the set $I_t$, there exist meromorphic mappings $c_{ti}\not\equiv 0\ (k_{t-1}+1\le i\le k_t)$ such that 
$$\sum_{i=k_{t-1}+1}^{k_t}c_{ti}\cdot (f,\tilde a_i)\in 
\biggl(\bigcup_{j=1}^{t-1}\bigcup_{i\in I_t}{(f,\tilde a_i)}\biggl)_{\mathcal {R}}.$$
Therefore, there exist meromorphic mappings $c_{ti}\in \mathcal {R}\ (1\le i\le k_{t-1})$ 
such that 
$$\sum_{i=1}^{k_t}c_{ti}\cdot (f,\tilde a_i)=0.$$
Define $c_{ti}=0$ for all $i>k_t.$ Then 
$$\sum_{i=1}^{k_l}c_{ti}\cdot (f,\tilde a_i)=0.$$
Since $\{c_{ti}(f,\tilde a_i)\}_{i=k_{t-1}+1}^{k_t}$ is $\mathcal {R}$-linearly independent, by again Lemma \ref{2.2} there exists an admissible set $\{\alpha_{t(k_{t-1}+1)},\ldots,\alpha_{tk_t}\}\subset \Z^m_+$ \ $(|\alpha_{ti}|\le k_t-k_{t-1}-1\le \rank_{\mathcal R}f-1=k)$ such that the matrix
$$\ A_t=\left (\mathcal {D}^{\alpha_{ti}}(c_{1j}(f,\tilde a_j));k_{t-1}+1\le i,j\le k_t \right)$$
has nonzero determinant.

Consider the following $(k_l-1)\times k_l$ matrix 
\begin{align*}
T&=\left (\mathcal {D}^{\alpha_{ti}}(c_{1j}(f,\tilde a_j));k_0+1\le i\le k_t,1\le j\le k_t \right)\\
\\
 &=\left [ \begin {array}{cccc}
\mathcal {D}^{\alpha_{12}}(c_{11}(f,\tilde a_1)) &\cdots & \mathcal {D}^{\alpha_{12}}(c_{1k_l}(f,\tilde a_{k_l})) \\
\mathcal {D}^{\alpha_{13}}(c_{11}(f,\tilde a_1)) &\cdots & \mathcal {D}^{\alpha_{13}}(c_{1k_l}(f,\tilde a_{k_l})) \\
\vdots &\vdots &\vdots\\
\mathcal {D}^{\alpha_{1k_1}}(c_{11}(f,\tilde a_1))&\cdots &\mathcal {D}^{\alpha_{1k_1}}(c_{1k_l}(f,\tilde a_{k_l}))\\
\mathcal {D}^{\alpha_{2k_1+1}}(c_{21}(f,\tilde a_1)) &\cdots & \mathcal {D}^{\alpha_{2k_1+1}}(c_{2k_l}(f,\tilde a_{k_l})) \\
\mathcal {D}^{\alpha_{2k_1+2}}(c_{21}(f,\tilde a_1)) &\cdots & \mathcal {D}^{\alpha_{2k_1+2}}(c_{2k_l}(f,\tilde a_{k_l})) \\
\vdots &\vdots &\vdots\\
\mathcal {D}^{\alpha_{2k_2}}(c_{21}(f,\tilde a_1))&\cdots &\mathcal {D}^{\alpha_{2k_2}}(c_{2k_t}(f,\tilde a_{k_l}))\\
\vdots &\vdots &\vdots\\
\mathcal {D}^{\alpha_{lk_{l-1}+1}}(c_{l1}(f,\tilde a_1)) &\cdots & \mathcal {D}^{\alpha_{lk_{l-1}+1}}(c_{lk_l}(f,\tilde a_{k_l})) \\
\mathcal {D}^{\alpha_{lk_{l-1}+2}}(c_{l1}(f,\tilde a_1)) &\cdots & \mathcal {D}^{\alpha_{lk_{l-1}+2}}(c_{lk_l}(f,\tilde a_{k_l})) \\
\vdots &\vdots &\vdots\\
\mathcal {D}^{\alpha_{lk_l}}(c_{lk}(f,\tilde a_1))&\cdots &\mathcal {D}^{\alpha_{lk_l}}(c_{lk_l}(f,\tilde a_{k_l}))\\
\end {array}
\right].
\end{align*}
\vskip0.3cm
Denote by $D_i$ the subsquare matrix obtained by deleting the $(i+1)$-th column 
of the minor  matrix  $T$. Since the sum of each row of $T$ is zero, we have
$$\det D_i={(-1)}^{i-1}\det D_1={(-1)}^{i-1}\prod_{j=1}^{l}\det A_j.$$

Since $\{a_i\}_{i=1}^q$ is in general position, we have 
$$\det (\tilde a_{ij}, \ 1\le i\le n+1,0\le j\le n )\not\equiv 0.$$
By solving the linear equation system $(f,\tilde a_i)=\tilde a_{i0}\cdot f_0+\ldots +\tilde a_{in}\cdot f_n \ (1\le i\le n+1),$
we obtain 
\begin{align}\label{+}f_v=\sum_{i=1}^{n+1}A_{vi}(f,\tilde a_{i})\ (A_{vi}\in\mathcal R)\text{ for each }0\le v \le n.
\end{align}
Put 
$\Psi(z)=\sum_{i=1}^{n+1}\sum_{v=0}^n |A_{vi}(z)|\ (z\in \C^m).$
Then
$$\ \ ||f(z)||\le \Psi (z)\cdot \max_{1\le i\le n+1}\bigl (|(f,\tilde a_i)(z)|\bigl )\le \Psi (z)\cdot \max_{1\le i\le q}\bigl (|(f,\tilde a_i)(z)|\bigl )\ (z\in \C^m),$$
and
\begin{align*}
\int\limits_{S(r)} \log^+\Psi (z) \eta &\le \sum_{i=1}^{n+1}\sum_{v=0}^n \int\limits_{S(r)}\log^+|A_{vi}(z)|\eta+O(1)\\
&\le \sum_{i=1}^{n+1}\sum_{v=0}^n T(r,A_{vi}) +O(1)\\
&=  O(\max_{1\le i\le q}T_{a_i}(r))+O(1).
\end{align*}

Fix $z_0 \in \C^m\setminus\bigcup_{j=1}^q\biggl (\supp (\nu^0_{(f,\tilde a_j)})\cup \supp (\nu^\infty_{(f,\tilde a_j)})\biggl ).$ Take $i\ (1\le i \le q)$ such that
$$|(f,\tilde a_i)(z_0)|=\max_{1\le j\le q}(|f,\tilde a_j)(z_0)|.$$
 Then
\begin{align*}
\dfrac{|\det D_1(z_0)|\cdot ||f(z_0)||}{\prod_{j=1}^{q}|(f,\tilde a_i)(z_0)|}
&=\dfrac{|\det D_i(z_0)|}{\prod_{\underset{j\ne i}{j=0}}^{q}|(f,\tilde a_j)(z_0)|}\cdot \biggl(\dfrac
{||f(z_0)||}{|(f,\tilde a_i)(z_0)|}\biggl)\\
&\le \Psi (z_0)\cdot \dfrac{|\det D_i(z_0)|}{\prod_{\underset{j\ne i}{j=1}}^{q}|(f,\tilde a_j)(z_0)|}.
\end{align*}
This implies that
\begin{align*}
\log\dfrac{|\det D_1(z_0)|.||f(z_0)||}{\prod_{j=1}^{q}|(f,\tilde a_j)(z_0)|}& \le\log^+\biggl ( \Psi (z_0)\cdot \biggl(\dfrac{|\det D_i(z_0)|}{\prod_{j=1,j\ne i}^{q}|(f,\tilde a_j)(z_0)|}\biggl)\biggl  )\\
& \le\log^+\biggl(\dfrac{|\det D_i(z_0)|}{\prod_{j=1,j\ne i}^{q}|(f,\tilde a_j)(z_0)|}\biggl)+\log^+\Psi (z_0).
\end{align*}

Thus, for each $z\in \C^m\setminus\bigcup_{j=1}^q\biggl (\supp (\nu^0_{(f,\tilde a_j)})\cup \supp (\nu^\infty_{(f,\tilde a_j)})\biggl ),$ we have
\begin{align*}
\log\dfrac{|\det D_1(z)|.||f(z)||}{\prod_{i=1}^{q}|(f,\tilde a_i)(z)|} \le\sum_{i=1}^{q}\log^+\biggl(\dfrac{|\det D_i(z)|}{\prod_{j=1,j\ne i}^{q}|(f,\tilde a_j)(z)|}\biggl)+\log^+ \Psi (z)
\end{align*}
Hence
\begin{align}\label{*}
\log ||f(z)||+\log \dfrac{|\det D_1(z)|}{\prod_{i=1}^{q}|(f,\tilde a_i)(z)|}\le \sum_{i=1}^{q}\log^+\biggl(\dfrac{|\det D_i(z)|}{\prod_{j=1,j\ne i}^{q}|(f,\tilde a_j)(z)|}\biggl)+\log^+ \Psi (z).
\end{align}
Note that   
\begin{align*}
\dfrac{\det D_i}{\prod_{j=1,j\ne i}^{q}(f,\tilde a_j)}
&=\dfrac{\det D_i/f_0^{q-1}}{\prod_{j=1,j\ne i}^{q}\biggl ((f,\tilde a_j)/f_0\biggl )}\\
&=\left [ \begin {array}{cccc}
\dfrac {\mathcal {D}^{\alpha_{12}}\biggl(\dfrac {c_{11}(f,\tilde a_1)}{f_0}\biggl)}{\dfrac{(f,\tilde a_1)}{f_0}} &\cdots &
\dfrac {\mathcal {D}^{\alpha_{12}}\biggl(\dfrac {c_{1k_l}(f,\tilde a_{k_l})}{f_0}\biggl)}{\dfrac{(f,\tilde a_{k_l})}{f_0}} \\
\vdots &\vdots &\vdots \\
\dfrac {\mathcal {D}^{\alpha_{lk_l}}\biggl(\dfrac {c_{l1}(f,\tilde a_1)}{f_0}\biggl)}{\dfrac{(f,\tilde a_1)}{f_0}} &\cdots &
\dfrac {\mathcal {D}^{\alpha_{lk_l}}\biggl(\dfrac {c_{lk_l}(f,\tilde a_{k_l})}{f_0}\biggl)}{\dfrac {(f,\tilde a_{k_l})}{f_0}}
\end {array}
\right]
\end{align*}

\quad (The determinant is counted after deleting the $i$-th column in the above matrix).

Each element of the above matrix has a form 
$$\dfrac {\mathcal {D}^{\alpha}\biggl(\dfrac {c(f,\tilde a_j)}{f_0}\biggl)}{\dfrac{(f,\tilde a_j)}{f_0}}=
\dfrac {\mathcal {D}^{\alpha}\biggl(\dfrac {c(f,\tilde a_j)}{f_0}\biggl)}{\dfrac{c(f,\tilde a_j)}{f_0}}\cdot c \ (c \in \mathcal {R}).$$
By lemma on logarithmic derivative lemma, we have
\begin{align*}
\biggl|\biggl| \quad\quad m \biggl(r,\dfrac {\mathcal {D}^{\alpha}\biggl(\dfrac {c(f,\tilde a_j)}{f_0}\biggl)}{\dfrac{(f,\tilde a_j)}{f_0}}\biggl)&\le
m \biggl(r,\dfrac {\mathcal {D}^{\alpha}\biggl(\dfrac {c(f,\tilde a_j)}{f_0}\biggl)}{\dfrac{c(f,\tilde a_j)}{f_0}}\biggl)+m(r,c)\\
&= O\biggl(\log^+T\biggl(r,\dfrac {c(f,\tilde a_j)}{f_0}\biggl)\biggl)+O(\max_{1\le i\le q}T(r,a_i))\\
&= O(\log^+T_f(r))+O(\max_{1\le i \le q}T(r,a_i)).
\end{align*}
This yields that
$$\biggl|\biggl| \quad m\left (r,\dfrac{\det D_i}{\prod_{j=1,j\ne i}^{q}(f,\tilde a_j)}\right )=
 O(\log^+T_f(r))+O(\max_{1\le j \le q}T_{a_j}(r))\ (1 \le i \le q).$$
Hence
$$\biggl|\biggl| \quad\quad \sum_{i=1}^{q} m\left (r,\dfrac{\det D_i}{\prod_{j=1,j\ne i}^{q}(f,\tilde a_j)}\right )= O(\log^+T_f(r))+O(\max_{1\le j \le q}T_{a_j}(r)).$$
 
Integrating both sides of the inequality (\ref{*}), we have
\begin{align*}
\biggl|\biggl|  \  \int_{S(r)}\log ||f|| \eta &+ \int_{S(r)}\log \biggl(\dfrac{|\det{D}_0|}{\prod_{i=1}^{q} |(f,\tilde a_i)|} \biggl)\eta\\
&\le \sum_{i=1}^{q} \int_{S(r)}\log^+ \biggl(\dfrac{|\det D_i|}{\prod_{j=1,j\ne i}^{q}
|(f,\tilde a_j)|}\biggl)\eta +\int_{S(r)}\log^+ \Psi(z)\eta\\
&= \sum_{i=1}^{q} m\biggl(r,\dfrac{\det D_i}{\prod_{j=1,j\ne i}^{q}(f,\tilde a_j)}\biggl)
+O(\max_{1\le i \le q}T_{a_i}(r))\\
&= O(\log^+T_f(r))+O(\max_{0\le i \le q-1}T_{a_i}(r)).
\end{align*}
Hence
$$||\ \ T_f(r)+ \int\limits_{S(r)} \text{log}\dfrac{|\det D_1|}{\prod_{i=1}^{q}|(f,\tilde a_i)|} \eta =O(\log^+T_f(r))+O(\max_{1\le i\le q}T_{a_i}(r)),\ \text {i.e, }$$
\begin{align}\nonumber
||\ T_f(r) &= \int\limits_{S(r)} \text{log}\dfrac{\prod_{i=1}^{q}|(f,\tilde a_i)|}{|\det D_1|} \eta+ O(\log^+T_f(r))+O(\max_{1\le i\le q}T_{a_i}(r))\\
\nonumber
&= \int\limits_{S(r)} \text{log}\prod_{i=1}^{q}|(f,\tilde a_i)|\eta- \int\limits_{S(r)} \text{log}|\det D_1| \eta +O(\log^+T_f(r))+O(\max_{1\le i\le q}T_{a_i}(r))\\
\label{3.3}
&\le N_{\prod_{i=1}^{q}(f, \tilde a_i)}(r)-N(r,\nu_{\det D_1})+ O(\log^+T_f(r))+O(\max_{1\le i\le q}T_{a_i}(r)).
\end{align}
\begin{claim}
$||\ N_{\prod_{i=1}^{q}(f, \tilde a_i)}(r)-N(r,\nu_{\det D_1})\le \sum_{i=1}^qN^{[k]}_{(f,a_i)}(r)+O(\max_{1\le i\le q}T_{a_i}(r)).$
\end{claim}
Indeed, fix $z\in\C^m\setminus I(f)$, where $I(f)=\{f_0=\cdots f_n=0\}$. We call $i_0$ the index satisfying 
$$\nu^0_{(f,\tilde a_{i_0})}(z)=\min_{1\le i\le n+1}\nu^0_{(f,\tilde a_i)}(z).$$
For each $i\ne i_0, i\in I_s$, we have
\begin{align*}
\nu^0_{\mathcal {D}^{\alpha_{sk_{s-1}+j}}(c_{si}(f,\tilde a_{i}))}(z)
&\ge\min_{\beta\in \Z_+^m \text { with } \alpha_{sk_{s-1}+j}-\beta \in\Z_+^m}
\{\nu^0_{\mathcal{D}^{\beta}c_{si}\mathcal D^{\alpha_{st_{s-1}+j}-\beta }(f,\tilde a_{i})}(z)\}\\
&\ge\min_{\beta\in \Z_+^n \text { with } \alpha_{sk_{s-1}+j}-\beta \in\Z_+^n}\bigl{\{}\max\{0,\nu^0_{(f,\tilde a_{i})}(z)-|\alpha_{sk_{s-1}+j}-\beta|\}\\
&\hspace{90pt}-(\beta+1)\nu^{\infty}_{c_{si}}(z)\bigl{\}}\\
&\ge\max\{0,\nu_{(f\tilde a_i)}^{0}(z)-k\}-(k+1)\nu_{c_{si}}^{\infty}(z)
\end{align*}
On the other hand, we also have
\begin{align*}
\nu^\infty_{\mathcal {D}^{\alpha_{sk_{s-1}+j}}(c_{si}(f,\tilde a_{i}))}(z)\le (|\alpha_{sk_{s-1}+j}|+1)\nu^\infty_{c_{si}(f,\tilde a_{i})}(z)\le (k+1)(\nu^\infty_{c_{si}}(z)+\nu^0_{a_{i0}}(z)).
\end{align*}
Thus
$$ \nu_{\mathcal {D}^{\alpha_{sk_{s-1}+j}}(c_{si}(f,\tilde a_{i}))}(z)\ge \max\{0,\nu_{(f\tilde a_i)}^{0}(z)-k\}-(k+1)\bigl (2\nu_{c_{si}}^{\infty}(z)+\nu^0_{a_{i0}}(z)\bigl )$$
Since each element of the matrix $D_{i_0}$ has a form $\mathcal {D}^{\alpha_{sk_{s-1}+j}}(c_{si}(f,\tilde a_{i}))\ (i\ne i_0)$, one estimates
\begin{align}\label{3.4}
\nu_{D_1}(z)=\nu_{D_{i_0}}(z)
\ge\sum_{i\ne i_0}\left (\max\{0,\nu_{(f\tilde a_i)}^{0}(z)-k\}-(k+1)\bigl (2\nu_{c_{si}}^{\infty}(z)+\nu^0_{a_{i0}}(z)\bigl )\right ).
\end{align}
We see that there exists $v_0\in\{0,\ldots,n\}$ with $f_{v_0}(z)\ne 0$. Then by (\ref{+}), there exists $i_1\in\{1,\ldots,n+1\}$ such that $A_{v_0i_1}(z)\cdot (f,\tilde a_{i_1})(z)\ne 0$. Thus
\begin{align}\label{3.5}
\nu^0_{(f,\tilde a_{i_0})}(z)\le \nu^0_{(f,\tilde a_{i_1})}(z)\le\nu^\infty_{A_{v_0i_1}}(z)\le\sum_{A_{vi}\not\equiv 0}\nu^\infty_{A_{vi}}(z).
\end{align}
Combining the inequalities (\ref{3.4}) and (\ref{3.5}), we have
\begin{align*}
\nu^{0}_{\prod_{i=1}^{q}(f,\tilde a_i)}(z)&-\nu_{\det D_1}(z)\\
&\le\sum_{i\ne i_0}\left (\min\{\nu_{(f,\tilde a_i)}^{0}(z),k\}+(k+1)\bigl (2\nu_{c_{si}}^{\infty}(z)+\nu^0_{a_{i0}}(z)\bigl )\right )+\sum_{A_{vi}\not\equiv 0}\nu^\infty_{A_{vi}}(z)\\
&\le\sum_{i=1}^q\left (\min\{\nu_{(f,\tilde a_i)}^{0}(z),k\}+(k+1)\bigl (2\nu_{c_{si}}^{\infty}(z)+\nu^0_{a_{i0}}(z)\bigl )\right )+\sum_{A_{vi}\not\equiv 0}\nu^\infty_{A_{vi}}(z),
\end{align*}
where the index $s$ of $c_{si}$ is taken so that $i\in I_s$. Integrating both sides of this inequality, we obtain
\begin{align}\nonumber
||\ \ N_{\prod_{i=1}^{q}(f, \tilde a_i)}(r)&-N(r,\nu_{\det D_1})\\
\nonumber
&\le \sum_{i=1}^q\left (N^{[k]}_{(f,\tilde a_i)}(r)+(k+1)\biggl (2N_{\frac{1}{c_{si}}}(r)+N_{a_{i0}}(r)\biggl )\right )+\sum_{A_{vi}\not\equiv 0}N_{{1}/{A_{vi}}}(r)\\
\label{3.6}
&= \sum_{i=1}^qN^{[k]}_{(f,a_i)}(r)+O(\max_{1\le i\le q}T_{a_i}(r)).
\end{align}
The claim is proved.

From the inequalities (\ref{3.3}) and the claim, we get
$$ ||\ \ T_f(r)\le \sum_{i=1}^qN^{[k]}_{(f,a_i)}(r)+O(\log^+T_f(r))+O(\max_{1\le i\le q}T_{a_i}(r)). $$
The lemma is proved.\hfill$\square$

\vskip0.2cm
\noindent
{\bf Proof of Theorem \ref{1.1}.}

(a). We denote by $\mathcal I$ the set of all permutations of $q-$tuple $(1,\ldots,q)$. For each element $I=(i_1,\ldots,i_q)\in\mathcal I$, we set
$$N_I=\{r\in\R^+;N^{[k]}_{(f,a_{i_1})}(r)\le\cdots\le N^{[k]}_{(f,a_{i_q})}(r)\}.$$

We now consider an element $I=(i_1,\ldots,i_q)$ of $\mathcal I$. We will construct subsets $I_t$ of the set $A_1=\{1,\ldots,{2n-k+2}\}$ as follows.

We choose a subset $I_1$ of $A$ which is the minimal subset of $A$ satisfying that $\{(f,\tilde a_{i_j})\}_{j\in I_1}$ is minimal over $\mathcal R$. If $\sharp I_1\ge n+1$ then we stop the process. 

Otherwise, set $A_2=A_1\setminus I_1$. We consider the following two cases:
\begin{itemize}
\item Case 1. Suppose that $\sharp A_2\ge n+1$. Since $\{\tilde a_{i_j}\}_{j\in A_2}$ is in general position, we have
$$ \left ((f,\tilde a_{i_j}); j\in A_2\right )_{\mathcal R}=\left (f_0,\ldots,f_n\right )_{\mathcal R}\supset \left ((f,\tilde a_{i_j}); j\in I_1\right )_{\mathcal R}\not\equiv 0.$$
\item  Case 2. Suppose that $\sharp A_2< n+1$. Then we have the following:
\begin{align*}
&\dim_{\mathcal R}\left ((f,\tilde a_{i_j}); j\in I_1\right )_{\mathcal R}\ge k+1-(n+1-\sharp I_1)=k-n+\sharp I_1,\\ 
& \dim_{\mathcal R}\left ((f,\tilde a_{i_j}); j\in A_2\right )_{\mathcal R}\ge k+1-(n+1-\sharp A_2)=k-n+\sharp A_2.
\end{align*}
We note that $\sharp I_1+\sharp A_2=2n-k+2$. Hence the above inequalities imply that
\begin{align*}
\dim_{\mathcal R}&\biggl (\bigl ((f,\tilde a_{i_j}); j\in I_1\bigl )_{\mathcal R}\cap\bigl ((f,\tilde a_{i_j}); j\in A_2\bigl )_{\mathcal R}\biggl )\\
&\ge\dim_{\mathcal R}\left ((f,\tilde a_{i_j}); j\in I_1\right )_{\mathcal R}+\dim_{\mathcal R}\left ((f,\tilde a_{i_j}); j\in A_2\right )_{\mathcal R}-(k+1)\\
&=k-n+\sharp I_1+k-n+\sharp A_2-(k+1)=1.
\end{align*}
\end{itemize}
Therefore, from the above two case, we see that
$$ \bigl ((f,\tilde a_{i_j}); j\in I_1\bigl )_{\mathcal R}\cap\bigl ((f,\tilde a_{i_j}); j\in A_2\bigl )_{\mathcal R}\ne \{0\}. $$
Therefore, we may chose a subset $I_2\subset A_2$ which is the minimal subset of $A_2$ satisfying that there exist nonzero meromorphic functions $c_i\in\mathcal R\ (i\in I_2)$,
$$\sum_{i\in I_2}c_i(f,\tilde a_i)\in \biggl(\bigcup_{i\in I_1}(f,\tilde a_i) \biggl)_{\mathcal {R}}.$$
By the minimality of the set $I_2$, the family $\{(f,\tilde a_{i_j})\}_{j\in I_2}$ is linearly independent over $\mathcal R$, and hence $\sharp I_2\le k+1$ and 
$$\sharp (I_2\cup I_2)\le\min\{2n-k+2, n+k+1\}.$$
If $\sharp (I_2\cup I_2)\ge n+1$ then we stop the process.

Otherwise, by repeating the above argument, we have a subset $I_3$ of $A_3=A_1\setminus (I_1\cup I_2)$, which satisfies the following:
\begin{itemize}
\item there exist nonzero meromorphic functions $c_i\in\mathcal R\ (i\in I_3)$ so that
$$\sum_{i\in I_3}c_i(f,\tilde a_i)\in \biggl(\bigcup_{i\in I_1\cup I_2}(f,\tilde a_i) \biggl)_{\mathcal {R}},$$
\item $\{(f,\tilde a_{i_j})\}_{j\in I_3}$ is linearly independent over $\mathcal R$,
\item $\sharp I_3\le k+1$ and  $\sharp (I_1\cup\cdots\cup I_3)\le \min\{2n-k+2, n+k+1\}$.
\end{itemize}

Continuing this process, we get the subsets $I_1,\ldots,I_l$, which satisfy: 
\begin{itemize}
\item $\{(f,\tilde a_{i_j})\}_{j\in I_1}$ is minimal over $\mathcal R$, $\{(f,\tilde a_{i_j})\}_{j\in I_t}$ is linearly independent over $\mathcal {R}\ (2\le t \le l), $
\item for any $2\le t\le l, j\in I_t,$ there exist meromorphic functions $c_j\in \mathcal {R}\setminus\{0\}$ such that 
$$\sum_{j\in I_t}c_j(f,\tilde a_{i_j})\in \biggl(\bigcup_{s=1}^{t-1}\bigcup_{j\in I_s}(f,\tilde a_{i_j}) \biggl)_{\mathcal {R}},$$
\item $n+1\le\sharp (I_1\cup\cdots\cup I_l)\le \min\{2n-k+2, n+k+1\}$.
\end{itemize}

Then the family of subsets $I_1,\ldots,I_t$ satisfies the assumptions of the Lemma \ref{3.1}. Therefore, we have
\begin{align*}
||\ T_f(r)\le\sum_{j\in J}N^{[k]}_{(f,a_{i_j})}+o(T_f(r))+O(\max_{1\le i \le q}T_{a_i}(r)),
\end{align*}
where $J=I_1\cup\cdots\cup I_l$. Then for all $r\in N_I$ (may be outside a finite Borel measure subset of $\R^+$) we have
\begin{align}\nonumber
||\ T_f(r)&\le\dfrac{\sharp J}{q-(2n-k+2)+\sharp J}\biggl (\sum_{j\in J}N^{[k]}_{(f,a_{i_j})}(r)+\sum_{j=2n-k+3}^qN^{[k]}_{(f,a_{i_j})}(r)\biggl )\\
\label{3.7}
&+o(T_f(r)) + O(\max_{1\le i \le q}T_{a_i}(r)).
\end{align}
Since $\sharp J\le 2n-k+2$, the above inequality implies that 
\begin{align}\label{3.8}
||\ T_f(r)\le\dfrac{2n-k+2}{q}\sum_{i=1}^qN^{[k]}_{(f,a_{i})}(r)+o(T_f(r))+O(\max_{1\le i \le q}T_{a_i}(r)),\quad r\in N_I.
\end{align}

We see that $\bigcup_{I\in\mathcal I}N_I=\R^+$ and the inequality (\ref{3.8}) holds for every $r\in N_I, I\in\mathcal I$. This yields that
$$ T_f(r)\le\dfrac{2n-k+2}{q}\sum_{i=1}^qN^{[k]}_{(f,a_{i})}(r)+o(T_f(r))+O(\max_{1\le i \le q}T_{a_i}(r)) $$ 
for all $r$ outside a finite Borel measure subset of $\R^+$. Thus
$$ ||\ \dfrac{q}{2n-k+2}T_f(r)\le\sum_{i=1}^qN^{[k]}_{(f,a_{i})}(r)+o(T_f(r))+O(\max_{1\le i \le q}T_{a_i}(r)). $$
The assertion (a) is proved. 

(b) We repeat the same argument as in the proof of the assertion (a). If $n+k+1>2n-k+1$ then the assertion (b) is a consequence  of the assertion (a). Then we now only consider the case where $n+k+1\le 2n-k+1$.

From (\ref{3.7}) with a note that $\sharp J\le n+k+2$, we have
\begin{align*}
||\ T_f(r)&\le\dfrac{n+k+1}{q-(2n-k+2)+n+k+1)}\sum_{i=1}^qN^{[k]}_{(f,a_{i})}(r)+o(T_f(r))+O(\max_{1\le i \le q}T_{a_i}(r))\\
&=\dfrac{n+k+1}{q-n+2k-1}\sum_{i=1}^qN^{[k]}_{(f,a_{i})}(r)+o(T_f(r))+O(\max_{1\le i \le q}T_{a_i}(r))\ r\in N_I.
\end{align*}
Repeating again the argument in the proof of assertion (a), we see that the above inequality holds for all $r\in\R^+$ outside a finite Borel measure set. Then the assertion (b) is proved. \hfill$\square$

\section{Proof of Theorem \ref{1.3}}

In order to prove Theorem \ref{1.3}, we need the following.

\noindent
{\bf 4.1.}\ Let $f:\C^m \to \P^n(\C)$ be a meromorphic mapping with a reduced representation $f=(f_0:\ldots :f_n)$.
Let $\{a_i\}_{i=1}^{q}$ be ``slowly'' (with respect to $f$) moving hyperplanes of $\P^n(\C)$ in general position such that 
$$\dim \{ z \in \C^m : (f,a_i)(z)=(f,a_j)(z)=0 \} \le m-2 \quad (1 \le i<j \le q).$$

For $M+1$ elements $f^0,\ldots ,f^M \in \mathcal {F}(f,\{a_j\}_{j=1}^q,1)$, we put
$$T(r)=\sum_{k=0}^M T(r,f^k).$$

Assume that $a_i$ has a reduced representation $a_i=(a_{i0}:\cdots :a_{in}).$ By changing the homogeneous coordinate system of $\P^n(\C),$ 
we may assume that $a_{i0}\not \equiv 0\ (1\le i \le q).$

We set $\ F^{jk}_{i}:=\dfrac{(f^k,a_j)}{(f^k,a_i)}\quad (1 \le i,j \le q, \ 0\le k\le M).$

\begin{lemma}\label{4.1}
Suppose that $q\ge 2n+1$. Then 
$$|| \ T_g(r)=O(T_f(r))\text { for each } g \in \mathcal {F}(f,\{a_i\}_{i=1}^q,1).$$
\end{lemma}
\noindent
{\bf Proof.} By Corollary \ref{1.2}(a), we have 
\begin{align*}
\parallel\ \dfrac {2q-n+1}{3(n+1)}T_g(r)&\le \sum_{i=1}^qN_{(g,a_i)}^{[n]}(r)+o(T_g(r)+T_f(r))\\
&\le n\sum_{i=1}^qN_{(g,a_i)}^{[1]}(r)+o(T_g(r)+T_f(r))\\
&= \sum_{i=1}^qnN_{(f,a_i)}^{[1]}(r)+o(T_g(r)+T_f(r))\\
&\le qn T_f(r)+o(T_g(r)+T_f(r)).
\end{align*}
Hence \quad $|| \quad T_g(r)=O(T_f(r)).$ \hfill$\square$

\begin{definition}[{see \cite[p. 138]{Fu3}}]
Let $F_0,\ldots ,F_M$ be nonzero meromorphic functions on $\C^m$, where $M\ge 1$. Take a set $\alpha:=(\alpha^0,\ldots ,\alpha^{M-1})$ whose components $\alpha^k$ are composed of $m$ nonnegative integers, and set $|\alpha|=|\alpha^0|+\ldots +|\alpha^{M-1}|.$ We define Cartan's auxiliary function by
$$\Phi^\alpha \equiv \Phi^\alpha(F_0,\ldots ,F_M):=F_0F_1\cdots F_M\left | \begin {array}{cccc}
1&1&\cdots &1\\
\mathcal {D}^{\alpha^0}(\frac {1}{F_0})&\mathcal {D}^{\alpha^0}(\frac {1}{F_1}) &\cdots &\mathcal {D}^{\alpha^0}(\frac {1}{F_M})\\
\vdots &\vdots &\vdots &\vdots \\
\mathcal {D}^{\alpha^{M-1}}(\frac {1}{F_0}) &\mathcal {D}^{\alpha^{M-1}}(\frac {1}{F_1}) &\cdots &\mathcal {D}^{\alpha^{M-1}}(\frac {1}{F_M})\\
\end {array}
\right|
$$
\end{definition}

\begin{lemma}[{see \cite[Proposition 3.4]{Fu3}}]\label{4.4}
 If $\Phi^\alpha(F,G,H)=0$ and $\Phi^\alpha(\frac {1}{F},\frac {1}{G},\frac {1}{H})=0$ for all $\alpha$ with $|\alpha|\le1$, then one of the following assertions holds :

(i) \ $F=G, G=H$ or $H=F$

(ii) \ $\frac {F}{G},\frac {G}{H}$ and $\frac {H}{F}$ are all constant.
\end{lemma}

\begin{lemma}[{see \cite[Lemma 4.7]{TQ05}}]\label{4.5}
Suppose that there exists $\Phi^\alpha=\Phi^\alpha(F_{i_0}^{j_00},\ldots , F_{i_0}^{j_0M})\not\equiv 0$ with $1\le i_0,j_0\le q,\  |\alpha|\le \dfrac {M(M-1)}{2}, \ d\ge |\alpha|.$ Assume that $\alpha$ is a minimal element such that 
$\Phi^\alpha(F_{i_0}^{j_00},\ldots , F_{i_0}^{j_0M})\not\equiv 0$. Then, for each $0 \le k \le M$, the following holds:
$$\parallel N_{(f^k,a_{j_0})}^{[d-|\alpha|]}(r)+M \sum_{j\ne{j_0,i_0}}N_{(f^k,a_j)}^{[1]}(r)\le N_{\Phi^\alpha}(r)
\le T(r)-M\cdot N^{[1]}_{(f^k,a_{i_0})}(r)+o(T(r)).$$
And hence
$$|| \quad N_{(f^k,a_{j_0})}^{[d-|\alpha|]}(r)+M \sum_{j\ne{j_0}}N_{(f^k,a_j)}^{[1]}(r)\le  T(r)+o(T(r)).$$
\end{lemma}

\vskip0.2cm
\noindent
{\bf 4.2. Proof of Theorem \ref{1.3}}

\noindent
a) Assume that $q>\frac{9n^2+9n+2}{2}$. Suppose that  there exist three distinct elements $f^0,f^1,f^2 \in  \mathcal {F}(f,\{a_j\}_{j=1}^{q},1).$

Suppose that there exist two indices $i,j\in\{1,\ldots ,q\}$ and  $\alpha=(\alpha_0,\alpha_1)\in (\Z_+^n)^2 $ with $|\alpha|\le 1$ such that $\Phi^\alpha (F_{j}^{i0},F_{j}^{i1},F_{j}^{i2}) \not \equiv 0$. By Lemma \ref{4.5}, we have 
$$2\sum_{t\ne i}N_{(f^0,a_t)}^{[1]}(r)\le T(r)+o(T_f(r)).$$
Hence, by Corollary \ref{1.2}(b) we have
\begin{align*}
\parallel T(r)&\ge \dfrac{2}{3}\sum_{k=1}^3 \sum_{t\ne i}N_{(f^k,a_t)}^{[1]}(r)+o(T_f(r))\ge \dfrac{2}{3n}\sum_{k=1}^3 \sum_{t\ne i}N_{(f^k,a_t)}^{[n]}(r)+o(T_f(r))\\
&\ge\dfrac{4(q-1)}{9n(n+1)}T(r)+o(T_f(r)).
\end{align*}
Letting $r\longrightarrow +\infty$, we get $1\ge\frac{4(q-1)}{9n(n+1)}$, i.e., $q\le\frac{9n^2+9n+4}{4}$. This is a contradiction.

Then for two indices $i,j$ $(1\le i<j\le q)$, we have 
$$\Phi^\alpha (F_{j}^{i0},F_{j}^{i1},F_{j}^{i2}) \equiv 0\text{ and }\Phi^\alpha (F_{i}^{j0},F_{i}^{j1},F_{i}^{j2}) \equiv 0$$
for all $\alpha=(\alpha_0,\alpha_1)\ \text { with }|\alpha|\le 1.$ By Lemma \ref{4.4}, there exists a constant $\lambda$ such that
$$F_{j}^{i0}=\lambda F_{j}^{i1},F_{j}^{i1}=\lambda F_{j}^{i2}, \text { or } F_{j}^{i2}=\lambda F_{j}^{i0}.$$
For instance, we assume that $F_{j}^{i0}=\lambda F_{j}^{i1}$. We will show that $ \lambda=1.$ 

Indeed, assume that $\lambda \ne 1$. Since $F_{j}^{i0}=F_{j}^{i1}$ on the set $\bigcup_{k\ne j}\{z : (f,a_k)(z)=0\},$ we have that 
$F_{j}^{i0}=F_{j}^{i1}=0$ on the set $\bigcup_{k\ne j}\{z : (f,a_k)(z)=0\}.$ Hence 
$\bigcup_{k\ne j}\{z : (f,a_k)(z)=0\} \subset \{z : (f,a_i)(z)=0\}.$ It follows that $\{z : (f,a_k)(z)=0\}=\emptyset \ (k\ne i,j).$ 
We obtain that
$$\parallel\dfrac{2(q-2)}{3(n+1)}T_f(r) \le \sum_{k\ne i,k\ne j}N_{(f,a_k)}^{[n]}(r)+o(T_f(r))=o(T_f(r)).$$
This is a contradiction. Thus $ \lambda =1\ (1\le i<j\le q).$

Define 
$$I_1= \{i\in \{1,\ldots ,q-1\}: F_{q}^{i0}=F_{q}^{i1}\},$$
$$I_2= \{i\in \{1,\ldots ,q-1\}: F_{q}^{i1}=F_{q}^{i2}\},$$
$$I_3= \{i\in \{1,\ldots ,q-1\}: F_{q}^{i2}=F_{q}^{i0}\}.$$
Since $\sharp (I_1\cup I_2\cup I_3)=\sharp \{1,\ldots ,q-1\}=q-1\ge 3n-2$, there exists $1\le k\le 3$ such that $\sharp \ I_k \ge n$.
Without loss of generality, we may assume that $\sharp \ I_1 \ge n$. This implies that $f^0=f^1$. This is a contradiction. 

Thus, we have $\sharp \  \mathcal {F}(f,\{a_i\}_{i=1}^{q},1)\le 2.$

\vskip0.2cm
\noindent
b)  Assume that $q>3n^2+n+2$. 

Take $g\in\mathcal{F}(f,\{a_i\}_{i=1}^{q},1).$ Suppose that $f\ne g.$ By changing indices if necessary, we may assume that
$$\underbrace{\dfrac{(f,a_1)}{(g,a_1)}\equiv \dfrac{(f,a_2)}{(g,a_2)}\equiv \cdots\equiv \dfrac{(f,a_{k_1})}
{(g,a_{k_1})}}_{\text { group } 1}\not\equiv
\underbrace{\dfrac{(f,a_{k_1+1})}{(g,a_{k_1+1})}\equiv \cdots\equiv \dfrac{(f,a_{k_2})}{(g,a_{k_2})}}_{\text { group } 2}$$
$$\not\equiv \underbrace{\dfrac{(f,a_{k_2+1})}{(g,a_{k_2+1})}\equiv \cdots\equiv \dfrac{(f,a_{k_3})}{(g,a_{k_3})}}_{\text { group } 3}\not\equiv \cdots\not\equiv \underbrace{\dfrac{(f,a_{k_{s-1}+1})}{(g,a_{k_{s-1}+1})}\equiv\cdots \equiv 
\dfrac{(f,a_{k_s})}{(g,a_{k_s})}}_{\text { group } s},$$
where $k_s=q.$ 

For each $1\le i \le q,$ we set
\begin{equation*}
\sigma (i)=
\begin{cases}
i+n& \text{ if $i+n\leq q$},\\
i+n-q&\text{ if  $i+n> q$}
\end{cases}
\end{equation*}
and  
$$P_i=(f,a_i)(g,a_{\sigma (i)})-(g,a_i)(f,a_{\sigma (i)}).$$
By supposition that $f\ne g$, the number of elements of each group is at most $n$. Hence $\dfrac{(f,a_i)}{(g,a_i)}$ and $\dfrac{(f,a_{\sigma (i)})}{(g,a_{\sigma (i)})}$ belong to distinct groups. This means that $P_i\not\equiv 0\ (1\le i\le q)$. 
 
Fix an index $i$ with $1\le i \le q.$ It is easy to see that
\begin{align*}
\nu_{P_i}(z) \ge \min\{\nu_{(f,a_i)},\nu_{(g,a_i)}\}+\min\{\nu_{(f,a_{\sigma (i)})},\nu_{(g,a_{\sigma (i)})}\}
+\sum_{\underset{v\ne i,\sigma (i)}{v=1}}^{q}\nu_{(f,a_v)}^{[1]}(z)
\end{align*}
outside a finite union of analytic sets of 
dimension $\le m-2.$ Since $\min\{a,b\}+n\ge\min\{a,n\}+\min\{b,n\}$ for all positive integers $a$ and $b$, 
the above inequality implies that
\begin{align*}
N_{P_i}(r)\geq \sum_{v=i,\sigma (i)}\left ( 
N^{[n]}_{(f,a_v)}(r)+N^{[n]}_{(g,a_v)}(r)-nN^{[1]}_{(f,a_v)}(r) \right )
+\sum_{\underset{v\ne i,\sigma (i)}{v=1}}^{q}N^{[1]}_{(f,a_v)}(r).
\end{align*} 
On the other hand, by the Jensen formula, we have
\begin{align*}
N_{P_i}(r)=&\int_{S(r)}\log |P_i|\eta + O(1)\\
\le &\int_{S(r)}\log (|(f,a_i)|^2+|(f,a_{\sigma (i)}|^2)^{\frac{1}{2}}\eta 
+ \int_{S(r)}\log (|(g,a_i)|^2+|(g,a_{\sigma (i)}|^2)^{\frac{1}{2}}\eta +O(1)\\
\le &T_f(r)+T_g(r) +o(T_f(r)).
\end{align*}
This implies that
\begin{align*}
T_f(r)+T_g(r)\ge &\sum_{v=i,\sigma (i)}\left ( 
N^{[n]}_{(f,a_v)}(r)+N^{[n]}_{(g,a_v)}(r)-nN^{[1]}_{(f,a_v)}(r) \right )\\
+ &\sum_{\underset{v\ne i,\sigma (i)}{v=1}}^{q}N^{[1]}_{(f,a_v)}(r)
+o(T_f(r)).
\end{align*}
Summing-up both sides of the above inequality over $i=1,\ldots ,q$ and by Corollary \ref{1.2}(b), we have
\begin{align*}
q(T_f(r)+T_g(r))\ge &2\sum_{v=i}^q\left ( 
N^{[n]}_{(f,a_v)}(r)+N^{[n]}_{(g,a_v)}(r) \right )\\
&+ (q-2n-2)\sum_{v=1}^{q}N^{[1]}_{(f,a_v)}(r)+o(T_f(r))\\
\ge & (2+\frac{q-2n-2}{2n})\sum_{v=i}^q\left ( 
N^{[n]}_{(f,a_v)}(r)+N^{[n]}_{(g,a_v)}(r) \right )+o(T_f(r))\\
\ge &(2+\frac{q-2n+2}{2n})\dfrac{2q}{3(n+1)}(T_f(r)+T_g(r))+o(T_f(r)).
\end{align*}
Letting $r \to \infty$, we get $q\ge (2+\frac{q-2n-2}{2n})\dfrac{2q}{3(n+1)}\Leftrightarrow q\le 3n^2+n+2.$ This is a contradiction.

Then $f=g$. This implies that $\sharp \mathcal{F}(f,\{a_i\}_{i=1}^{q},1)=1$. The theorem is proved. \hfill $\square$

\end{document}